\documentclass[10pt,twoside]{siamltex}
\usepackage{amsfonts,epsfig}

\def\cvd{~\vbox{\hrule\hbox{%
     \vrule height1.3ex\hskip0.8ex\vrule}\hrule } }

\newtheorem{formula}[theorem]{Formula}

\newcommand{\fq}{\mathbb{F}_q}
\newcommand{\pr}{\mathbf{P}} 
\newcommand{\E}{\mathbf{E}}
\newcommand{\wt}{\mathrm{wt}\,}
\newcommand{\rk}{\mathrm{rk}\,}
\newcommand{\ra}{\rightarrow}
\newcommand{\var}{\mathrm{var}\,}

\begin{document}
   
\bibliographystyle{plain}
\title{
Weight and Rank of Matrices over Finite Fields
}

\author{
Theresa Migler \thanks{Department of Mathematics, California Polytechnic State
University, San Luis Obispo, California 93407, USA
(tmigler@calpoly.edu, kmorriso@calpoly.edu, mogle@calpoly.edu).}
\and
Kent E. Morrison\footnotemark[2]
\and
Mitchell Ogle\footnotemark[2]
}
\pagestyle{myheadings}
\markboth{T. \ Migler, K.E.\ Morrison, M.\ Ogle}{Weight and Rank}
\maketitle

\begin{abstract}
   Define the weight of a matrix to be the number of non-zero entries. 
   One would like to count $m$ by $n$ matrices over a finite field by
   their weight and rank.  This is equivalent to determining the
   probability distribution of the weight while conditioning on the
   rank.  The complete answer to this question is far from finished. 
   As a step in that direction this paper finds a closed form for the
   average weight of an $m$ by $n$ matrix of rank $k$ over the finite
   field with $q$ elements.  The formula is a simple algebraic
   expression in $m$, $n$, $k$, and $q$.  For rank one matrices a
   complete description of the weight distribution is given and a
   central limit theorem is proved.
\end{abstract}

\section{Introduction}
For an $m \times n$ matrix $A$ over the finite field $\fq$ the
\textbf{weight} of $A$, denoted $\wt A$, is the number of
non-zero entries.  In the Hamming metric on matrices it is the 
distance from $A$ to $0$.

There is some relationship between the rank and the weight of a
matrix.  For example, if $\wt A =1$, then $\rk A =1$, and if $\rk A
=k$, then $\wt A \geq k$.  On the other hand, there are matrices of
rank one and maximal weight $mn$, such as a matrix with every entry a 
one.  In this article we determine the
average weight of rank $k$ matrices in terms of $k$, $m$, $n$, and
$q$.  Without fixing the rank, the average weight of $m \times n$
matrices is $mn(1-\frac{1}{q})$ and the weight has a binomial
distribution.  However, the full probability distribution of the
weight for matrices of rank $k$ is yet to be determined.

The tools are those of elementary combinatorics and linear algebra. 
Nothing special is used from the theory of finite fields other than 
the understanding that the fundamental ideas of linear algebra work 
over all fields and not just the real or complex numbers.

We need a modest amount of background material. We use  
three basic formulas. 
\begin{formula}
The number of ordered $k$-tuples of linearly 
independent vectors in ${\fq}^{n}$is
\[    (q^{n}-1)(q^{n}-q)(q^{n}-q^{2})\cdots(q^{n}-q^{k-1})
\]    
\end{formula}
\proof  The first vector is any non-zero vector and each succeeding 
vector must avoid the span of the previous vectors.  \cvd

\begin{formula}
    The number of $k$-dimensional subspaces of ${\fq}^{n}$ is given by 
the $q$-binomial coefficient
\[
{n \brack{k}}_q = \frac{(q^{n} - 1)(q^{n} - q)(q^{n} - q^2) \cdots (q^{n} - q^{k-1}) }
        {(q^{k} - 1)(q^{k} - q)(q^{k} - q^2) \cdots (q^{k} - q^{k-1})}
\]
\end{formula}
\proof  The numerator is the number of bases of all $k$-dimensional 
subspaces, while the denominator is the number of bases of any given 
subspace.  \cvd
\begin{formula}
    The number of $m \times n$ matrices of rank $k$ is
    \begin{eqnarray*}
	& & {m \brack{k}}_q (q^{n}-1)(q^{n}-q)\cdots(q^{n}-q^{k-1})\\
	& = &{n \brack{k}}_q (q^{m}-1)(q^{m}-q)\cdots(q^{m}-q^{k-1}) \\
	& = & \frac{(q^{m}-1)(q^{m}-q)\cdots(q^{m}-q^{k-1})(q^{n}-1)
	            (q^{n}-q)\cdots(q^{n}-q^{k-1})}
		    {(q^{k} - 1)(q^{k} - q)(q^{k} - q^2) \cdots (q^{k} - q^{k-1})}
    \end{eqnarray*}
\end{formula}
\proof  For a fixed $k$-dimensional subspace $W \subset {\fq}^{m}$, the 
number of matrices with $W$ as the column space is equal to the 
number of $k \times n$ matrices of rank $k$. Such a matrix is given by 
the $k$ linearly independent row vectors of length $n$. The number of 
those is given by Formula 1. The number of $k$-dimensional subspaces 
of ${\fq}^{m}$ is ${m \brack{k}}_q$ and the product is the number of 
rank $k$ matrices given in the first line. By the same reasoning, the 
second line counts the number of $n \times m$ matrices of rank $k$, 
which is the same.
\cvd  

A special case of Formula 3 is worth noting. The number of 
invertible $n \times n$ matrices is
\[
    (q^{n}-1)(q^{n}-q)\cdots(q^{n}-q^{n-1})
\]

\section{Average Weight}
The average weight of a rank $k$ matrix is the sum of the average 
weights of the entries, and the average weight of the $ij$ entry is 
the probability that the entry is not zero:
\[
   \E(\mathrm{wt} \, A)=\sum_{i,j}\pr(a_{ij} \neq 0)
\]
\begin{theorem}
The probability that $a_{ij} \neq 0$ for a rank $k$ matrix 
$A$ is the same for all $i$ and $j$.
\end{theorem}
\proof The probability that the $ij$ entry is not zero is the quotient
whose numerator is the number of matrices $A$ of rank $k$ with $a_{ij}\neq 
0$, and whose denominator is the number of matrices of rank $k$.
Consider the map on the $m \times
n$ matrices that switches rows 1 and $i$ and switches columns 1 and
$j$.  This map preserves rank and gives a bijection between the subset
of matrices of rank $k$ with a non-zero in the 1,1 location and the subset
of matrices of rank $k$ with a non-zero in the $i,j$ location.  Thus,
$\pr(a_{ij} \neq 0)=\pr[a_{11} \neq 0)$.  \cvd 

Call this common value \textbf{the average
weight per entry}.  Now we focus on the upper left corner of the
matrices of rank $k$.  Our analysis depends on the reduced row echelon
form. We recall the definition \cite{Lay97}.
\begin{definition} {\rm
   A rectangular matrix is in  \emph{row echelon form} if it 
   has the following three properties:
   \begin{enumerate}
      \item  All non-zero rows are above any rows of all zeros.
   
      \item  Each leading entry of a row is in a column to the right of 
the leading entry of the row above it.
   
      \item  All entries in a column below a leading entry are zero.
   \end{enumerate}
      If a matrix in echelon form satisfies the following additional conditions, 
then it is in  \emph{reduced row echelon form}:
   \begin{enumerate}
      \setcounter{enumi}{3}
      \item  The leading entry in each nonzero row is 1
   
      \item  Each leading 1 is the only non-zero entry in its column.
   \end{enumerate}
}
\end{definition}

The $k \times n$ matrices in reduced row echelon form correspond
bijectively with the $k$-dimensional subspaces of ${\fq}^{n}$.  The rows
of the matrix give a basis of the corresponding subspace.  When an $m
\times n$ matrix $A$ is reduced to reduced row echelon form by row
operations, the result is an $m \times n$ matrix whose first $k$ rows
form a basis of the row space of $A$.  Let $R$ be the $k \times n$
matrix consisting of the $k$ non-zero rows of the reduced form.  Then
$A$ and all matrices with the same row space can be constructed from
$R$ by multiplying $R$ on the left by an $m \times k$ matrix $C$ of 
rank $k$. The matrix $C$ is unique. This gives a factorization of $A$ 
as $A=CR$. In terms of the associated linear maps, $A$ is a linear map 
from ${\fq}^{n}$ to ${\fq}^{m}$, which factors into a surjective map 
onto ${\fq}^{k}$ followed by an injective map from ${\fq}^{k}$ to ${\fq}^{m}$. 
Recall that knowing the row space of a matrix is equivalent to knowing 
the kernel of the associated linear map. Thus, when the reduced 
matrix $R$ is held fixed and $C$ is varied, the product $CR$ gives 
all maps with the same row space (i.e. same kernel).

\begin{theorem}
   For $m \times n$ matrices of rank $k$, the average weight per entry 
   is \[  \frac{\left(1-\frac{1}{q} \right)\left(1-\frac{1}{q^{k}} \right)}
     {\left(1-\frac{1}{q^{m}} \right)\left(1-\frac{1}{q^{n}} \right)} 
     \] 
\end{theorem}
\proof  We consider the random selection of a rank $k$ matrix such 
that each such matrix is equally probable. With the factorization 
$A=CR$, this can be done by selecting $C$ uniformly from all $m 
\times k$ matrices of rank $k$ and by selecting $R$ independently from 
among all reduced row echelon matrices, which is the same as selecting 
the row space uniformly from all $k$-dimensional subspaces of ${\fq}^{n}$.
The upper left corner of $A$ is
$a_{11}=c_{11}r_{11}+c_{12}r_{21}+\cdots+c_{1k}r_{k1}$.  But because
$R$ is in reduced row echelon form, the first column of $R$ is either
all zeros or has a leading 1 followed by zeros.  Thus,
$a_{11}=c_{11}r_{11}$.  In order for $a_{11}$ to be non-zero, both
$r_{11}$ and $c_{11}$ must be non-zero.  Since the selection of $R$ is 
independent of the selection of $C$,
\[  
   \pr(a_{11} \neq 0) = \pr(c_{11} \neq 0)\pr(r_{11} \neq 0)
\]
The columns of $C$ are $k$ linearly independent vectors of length $m$ 
and so the first column is not the zero vector. That means there are 
$q^{m}-1$ possible first column vectors. There are $(q-1)$ choices 
for $c_{11} \neq 0$ and $q^{m-1}$ choices for the remaining entries of 
the first column. Therefore,
\[ 
   \pr(c_{11} \neq 0) = \frac{(q-1)q^{m-1}}{q^{m}-1}=
        \frac{q^{m}-q^{m-1}}{q^{m}-1}
\]
Now $r_{11}=0 \mbox{ or } 1$, and $r_{11}=0$ when the row space of $R$
contains nothing in the direction of the vector $(1,0,0,\ldots,0)$,
which is to say that the row space is contained in the
$(n-1)$-dimensional space $\{(0,x_{2},\ldots,x_{n}) \}$.  Therefore,
\[
    \pr(r_{11}=0) =\frac{{n-1 \brack k}_{q}}{{n \brack k}_{q}}
\]
\[ 
   \pr(r_{11} \neq 0) =1-\frac{{n-1 \brack k}_{q}}{{n \brack k}_{q}}
\]
Using Formula 2 one easily obtains
\[ 
   \frac{{n-1 \brack k}_{q}}{{n \brack k}_{q}} =
    \frac{q^{n-k}-1}{q^{n}-1}
\]
Putting these results together we have
\begin{eqnarray*}  
   \pr(a_{11} \neq 0) & = &  \frac{q^{m}-q^{m-1}}{q^{m}-1}
                 \left( 1-  \frac{q^{n-k}-1}{q^{n}-1}  \right) \\
     &=& \frac{(q^{m}-q^{m-1})(q^{n}-q^{n-k})}{(q^{m}-1)(q^{n}-1)} \\
     &=& \frac{\left(1-\frac{1}{q} \right)\left(1-\frac{1}{q^{k}} \right)}
     {\left(1-\frac{1}{q^{m}} \right)\left(1-\frac{1}{q^{n}} \right)} 
     \mbox{\qquad \qquad}  \cvd
\end{eqnarray*}  

With this result we have a clear picture of the effect of the
parameters $k$, $m$, and $n$ on the average weight.  The factor $1 -
1/q$ is the average weight per entry without the rank condition, in
which case the matrix size does not matter.  Note that with $m$ and
$n$ fixed, it is more likely for an entry to be non-zero for matrices
of larger rank, an intuitively plausible result because both weight
and rank are some measure of distance from the zero matrix.  Also, one
can see that as $m$, $n$, and $k$ simultaneously go to infinity, the
probability approaches $1 - 1/q$, which is again the unconditioned
probability.  For invertible matrices of size $n$ (i.e. $k=m=n$)
the probability of a non-zero entry is
\[ 
   \frac{1-\frac{1}{q}}{1-\frac{1}{q^{n}}}
\]

\section{Weight of Rank One Matrices}
For the matrices of rank one a more complete analysis of the weight
distribution is possible.  In this case $C$ is a non-zero column
vector of length $m$ and $R$ is a non-zero row vector of length $n$
whose leading non-zero entry is 1.  The rank one matrix $A=CR$ is
given by $a_{ij}=c_{i}r_{j}$, and so the weight of $A$ is the product
of the weights of $C$ and $R$. The weight of $C$ has a binomial 
distribution conditioned on the weight being positive (the entries of 
$C$ cannot all be zero)
\[  
  \pr(\mathrm{wt} \, C=\mu)=\frac{{m \choose \mu}(q-1)^{\mu}q^{m-\mu}}{(q^{m}-1)}
\]
Likewise for $R$ the weight distribution is given by
\[  
  \pr(\mathrm{wt} \, R=\nu)=\frac{{n \choose \nu}(q-1)^{\nu}q^{n-\nu}}{(q^{n}-1)}
\]
(To select a random $R$, choose a random non-zero vector of length $n$ and 
then scale it to make the leading non-zero entry 1. The scaling does 
not change the weight.)

The weight on rank one matrices is the product of these two binomial 
random variables, each conditioned to be positive.
\begin{eqnarray} 
   \pr(\mathrm{wt} \, A=\omega) &=& \sum_{\mu \nu =\omega}\pr(\mathrm{wt} 
   \, C=\mu)\pr(\mathrm{wt} \, R=\nu)\\
   &=& \sum_{\mu \nu =\omega}{m \choose \mu}{n \choose \nu}
      \frac{(q-1)^{m+n-\mu-\nu}q^{\mu+\nu}}{(q^{m}-1)(q^{n}-1)}
\end{eqnarray}
Not all weights between 1 and $mn$ occur for rank 1 matrices since the 
weight is a product with one factor no greater than $m$ and the other 
factor no greater than $n$. Plots of actual probability distributions
show spikes and gaps. Plots of cumulative distributions are smoother 
and lead us to expect a limiting normal distribution. See Figures 1 
and 2.

\begin{figure}[h]  
\centering
\epsfig{file=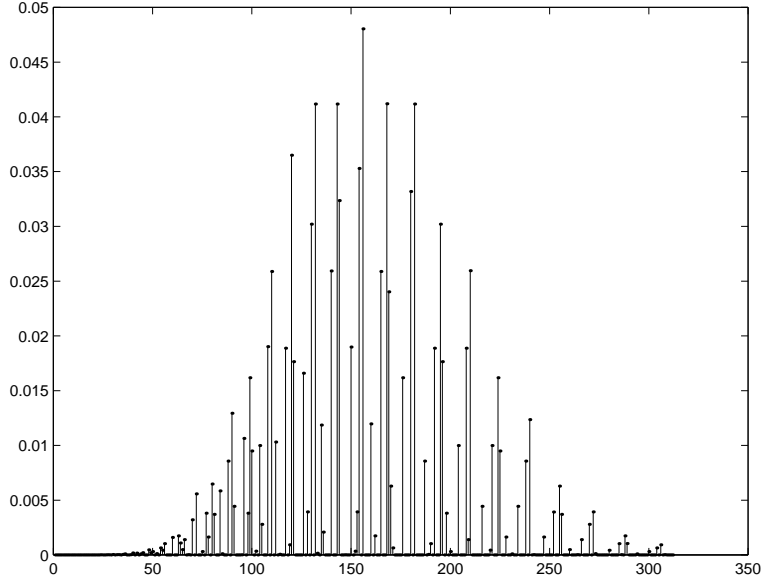,width=4in}
\caption{Distribution for the weight of rank 1 
matrices, $m=n=25$, $q=2$.  }
\end{figure}

\begin{figure}[h]  
\centering
\epsfig{file=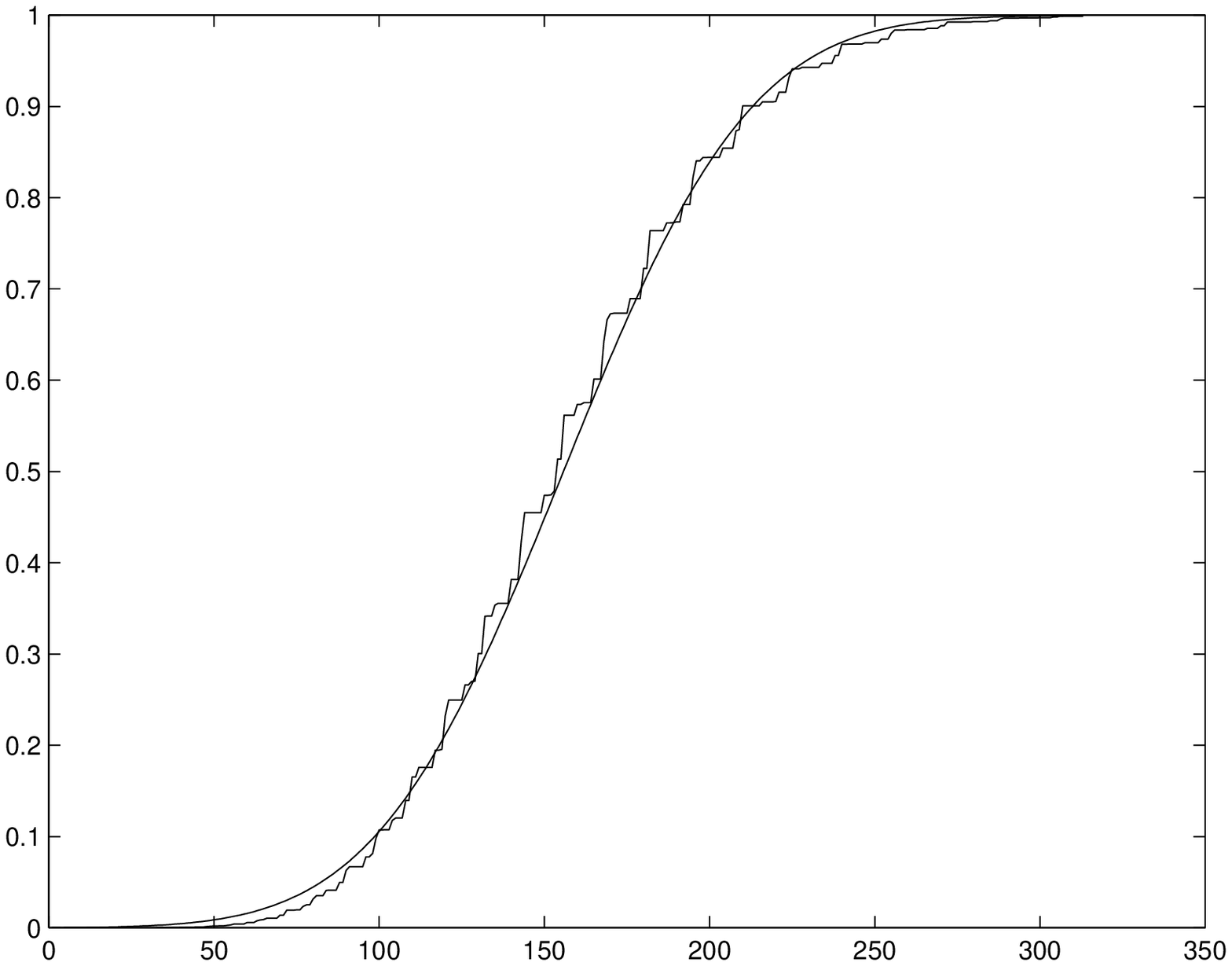,width=4in}
\caption{Cumulative frequency distribution for the weight of rank 1 
matrices, $m=n=25$, $q=2$. The smooth curve is the normal cdf with the same 
mean ($\approx 156.25$) and standard deviation ($\approx 44.63$). }
\end{figure}

\begin{theorem} As $m$ or $n$ goes to infinity, the weight 
distribution of rank one matrices approaches a normal distribution.
\end{theorem}
\proof The weight random variable for rank one matrices of size $m
\times n$ is the product of independent binomial random variables
conditioned on being positive.  Define $W=XY$, where $X=\sum_{1 \leq i
\leq m}X_{i}$, $Y=\sum_{1 \leq j \leq n}Y_{j}$, and $X_{i}$ and
$Y_{j}$ are independent Bernoulli random variables with probability
$1/q$ of being 0.  Then $W$ is the sum of $m$ independent identically
distributed random variables $X_{i}Y$. Conditioning $W$ on $W>0$ is
the weight of rank one matrices.  By the Central Limit Theorem the
distribution of $W$ converges, as $m \ra \infty$, to a normal
distribution after suitable scaling.  Now conditioning on $W$ being
positive does not change this result because the probability that
$W>0$ is $1-q^{-m}$, which goes to 1 as $m \ra \infty$.  \cvd

Now to compute the mean and variance of the weight, let $W=XY$ as
before but without conditioning on $X$ or $Y$ being positive.  Then
$\E(W) = mn(1-1/q)^{2}$ and

\[ \E(W^{2})=\E(X^{2}Y^{2}) = \E\left( \left(\sum_{i}X_{i}\right)^{2}
  \left(\sum_{j}Y_{j}\right)^{2} \right) \] Expanding and using the
  independence of the random variables $X_{i}, Y_{j}$ and the fact that 
  $X_{i}^{2}=X_{i}$ and $Y_{j}^{2}=Y_{j}$, we get
 \begin{eqnarray*} \E(W^{2}) &=& mn\left(1-\frac{1}{q}\right)^{2}
                     +mn(m+n-2)\left(1-\frac{1}{q}\right)^{3} \\
	& & \mbox{\quad}  + m(m-1)n(n-1)\left(1-\frac{1}{q}\right)^{4} 
 \end{eqnarray*}
The variance of the weight is
\begin{eqnarray*}
   \var(W|W>0) &=& \E(W^{2}|W>0) - {\E(W|W>0)}^{2} \\
      &=& \frac{\E(W^{2})}{\pr(W>0)}-\left( \frac{\E(W)}{\pr(W>0)}\right)^{2}
\end{eqnarray*} 
Furthermore 
\[ \pr(W>0) =  \pr(X > 0) \pr(Y>0) = (1-\frac{1}{q^{m}})(1-\frac{1}{q^{n}}) \]
Combining these expressions we get
\begin{eqnarray*}
   \lefteqn{\var(W | W>0) = } \\
  & & \frac{mn\left(1-\frac{1}{q}\right)^{2}
                     +mn(m+n-2)\left(1-\frac{1}{q}\right)^{3}
		     +m(m-1)n(n-1)\left(1-\frac{1}{q}\right)^{4}}
     {\left(1-\frac{1}{q^{m}}\right)\left(1-\frac{1}{q^{n}}\right)}  \\
  & & \mbox{\quad} - \frac{\left(mn\left(1-\frac{1}{q} \right) \right)^{2}}
    {\left(1-\frac{1}{q^{m}} \right)^{2}\left(1-\frac{1}{q^{n}} \right)^{2}}
\end{eqnarray*}
We get a good approximation to the variance of the weight when $m$ and 
$n$ using the unconditioned $W$, essentially the factors in the 
denominator by 1. Thus,
\begin{eqnarray*} \var(W) &= &\E(W^{2})-\E(W)^{2} \\
   &=&  mn\left(1-\frac{1}{q}\right)^{2}
        +mn(m+n-2)\left(1-\frac{1}{q}\right)^{3}  \\
   & & \mbox{\quad} + m(m-1)n(n-1)\left(1-\frac{1}{q}\right)^{4}  
        -\left(mn\left(1-\frac{1}{q}\right)\right)^{2} \\
\end{eqnarray*}
which can be simplified to give
\begin{eqnarray*}
   \var(W) &=& 
   mn(1-n-m)\left(1-\frac{1}{q}\right)^{4} \\
   & & \mbox{\quad} + mn(n+m-2)\left(1-\frac{1}{q}\right)^{3}  \\
   & &\mbox{\quad} + mn\left(1-\frac{1}{q}\right)^{2}
\end{eqnarray*}

From this we can see, for example, that for $m \approx n$, $m, n \ra \infty$, 
the variance is on the order of $n^{2}$ and the standard deviation is 
of order $n$.


\begin{thebibliography}{9}
    \bibitem{Lay97}  D. C. Lay.
      \newblock \emph{Linear Algebra and Its Applications}, second 
      edition.
      \newblock Addison-Wesley, Reading, Massachusetts, 1997.
\end{thebibliography}
\end{document}